\section{Background}\label{sec:background}
\subsection{Elliptic curves}
 The degree of a polynomial $A\in\mathbb{F}_q[t]$ is denoted by deg$(A)$. We define $N_E$, the conductor of $E$, by
\[N_E = \prod_{\substack{P\in\mathbb{F}_q[t] \\ P \text{ prime}}}P^{f_P(E)}\]
where
\begin{align*}
f_P(E):=\begin{cases}
 0, & \text{if } E \text{ has good reduction at } P \\
 1, &  \text{if } E \text{ has multiplicative reduction at } P  \\
 2, &  \text{if } E \text{ has additive reduction at } P
  \end{cases} 
\end{align*}
including the prime at infinity. A good reduction means that $E/P$, the reduction of $E$ at a prime $P\in\mathbb{F}_q[t]$, is an elliptic curve over $\mathbb{F}_q[t]/(P)\cong \mathbb{F}_{q^{\text{deg}(P)}}$, which is a finite field. By the Weil's conjectures, its zeta function is equal to
\[ \mathcal{Z}(E/P,u)=\frac{(1-\alpha_P(E)u)(1-\overline{\alpha_P}(E)u)}{(1-u)(1-q^{\text{deg}(P)}u)}\]
where $\alpha_P(E)$ and $\overline{\alpha_P}(E)$ may be swapped. We call those values the Frobenius eigenvalues of $E$ at $P$ and since we fixed $E$, we simply write them as $\alpha_P$ and $\overline{\alpha}_P$. For primes of good reduction, they may be computed using
\[ a_P:= \alpha_P + \overline{\alpha}_P = q^{\text{deg}(P)}+1-\#(E/P)\]
where $\#(E/P)$ is the number of $\mathbb{F}_{q^{\text{deg}(P)}}$-points of $E/P$. When $E/P$ isn't an elliptic curve, we say that the reduction is bad. If the reduction is multiplicative, we have $a_P=\pm 1$, and if the reduction is additive, we have $a_P=0$. We will often use the fact that $|\alpha_P|=q^{\text{deg}(P)/2}$ from the Weil's conjectures for primes of good reduction.
\subsection{$L$-functions}
Let $\chi$ be a primitive Dirichlet character on $\mathbb{F}_q[t]$ (see \cite{rosen:2013} for background). Its $L$-function is defined by
\[ L(\chi,u):= \prod_{P} (1-\chi(P)  u^{\text{deg}(P)})^{-1}\]
where the product is over the primes of $\mathbb{F}_q[t]$ excluding the prime at infinity. \\ \\
The $L$-function of E is defined by
\[ L(E,u):= \prod_{P\nmid N_E} (1-\alpha_P  u^{\text{deg}(P)})^{-1}(1-\overline{\alpha}_P  u^{\text{deg}(P)})^{-1} \prod_{P\mid N_E} (1-a_P u^{\text{deg}(P)})^{-1}\]
and the $L$-function of $E$ twisted by $\chi$ is defined by
\[ L(E\otimes \chi,u):= \prod_{P\nmid N_E} (1-\chi(P)\alpha_P  u^{\text{deg}(P)})^{-1}(1-\chi(P)\overline{\alpha}_P  u^{\text{deg}(P)})^{-1} \prod_{P\mid N_E} (1-\chi(P)a_P u^{\text{deg}(P)})^{-1}\]
excluding the prime at infinity.\\ \\
A Dirichlet character is uniquely factored into
\begin{equation}\label{eq:dc}
 \chi=\hat{\chi}\chi_0
 \end{equation}
where $\hat{\chi}$ is primitive and $\chi_0$ is principal with minimal modulus. For example, if $\chi$ is primitive, then $\chi_0\equiv 1$. If $\chi$ is principal, then $\hat{\chi} \equiv 1$. A character is said to be even if it is trivial on $\mathbb{F}_q$. It is called odd otherwise.\\ \\
We also split $\mathcal{H}_N^*=\mathcal{H}_{N,+}\cup \mathcal{H}_{N,-}$ depending on if the sign of the functional equation of $L(E\otimes\chi_D,u)$ for $D\in\mathcal{H}_N^*$ is 1 or $-1$ respectively.
\subsection{Gauss sums}
As done by Hayes \cite{hayes:1966}, we define the following function on $\mathbb{F}_q(t)$
\[ e_q(A) := e^{\frac{2\pi i \text{tr}_{\mathbb{F}_q/\mathbb{F}_p}(A_1)}{p}} \]
where $A_1$ is the coefficient of $1/t$ in the Laurent expansion of $A$. \\ The Gauss sum of a primitive character $\chi$ of conductor $F$  is defined as
\[ G(\chi) := \sum_{A\,\text{mod}\,F} \chi(A) e_q\left(\frac{A}{F}\right) \]
and it does not depend on the choice of representatives. \\
We also define
\[ \tau(\chi):=\sum_{a \in \mathbb{F}_q^*} \chi(a)e^{2\pi i \text{tr}_{\mathbb{F}_q/\mathbb{F}_p}(a)/p}.\]
We denote by $\omega(\chi)$ the sign of the functional equation of $L(\chi,u)$ and it is given in Lemma \ref{lm:duality}. When $\chi$ is quadratic, $\omega(\chi)$ is always 1 (see \cite{Rudnick:2008}). \\ \\
We denote the orthogonal group of dimension $M$  by $O(M)$. The subgroup of matrices with determinant 1 is denoted by $SO(M)$. Two different symmetry types exist depending on the parity of $M$. They are denoted by $SO(\text{even})$ and $SO(\text{odd})$.

\section{Duality}
\begin{proposition}[Riemann Hypothesis]\label{prop:rh}
Let $\chi$ be a primitive Dirichlet character of conductor $F\neq1$ on $\mathbb{F}_q[t]$. Then
\[  L(\chi,u) = (1-u)^\lambda \prod_{j=1}^M (1-q^{1/2}e^{i\theta_j}u) \]
where $M=\emph{deg}(F)-1$ and $\lambda=0$ if $\chi$ is odd or else $M=\emph{deg}(F)-2$ and $\lambda=1$.
\end{proposition}
\begin{proof}
We have by Theorem 9.16B of \cite{rosen:2013}
\[ L(\chi,u)K_\infty(\chi,u)= \prod_{j=1}^M (1-q^{1/2}e^{i\theta_j}u)\]
where
\[ K_\infty(\chi,u)=(1-\chi(P_\infty)u)^{-1}\]
is the Euler factor of the prime at infinity. This implies $L(\chi,u)$ have at most one zero on $|u|=1$ and all others are on $|u|=q^{-1/2}$. \\ \\
We have the expansion 
\[ L(\chi,u)=\sum_{h=0}^\infty\left(\sum_{D \in \mathcal{M}_h}\chi(D)\right)u^h.\]
When $h\geq \text{deg}(F)$, the polynomials become equidistributed modulo $F$ so the coefficient of $u^h$ is zero. We also have by \cite{rosen:2013} that 
\[ \sum_{D \in \mathcal{M}_{\text{deg}(F)-1}}\chi(D) \neq 0\]
meaning that $L(\chi,u)$ is a polynomial of degree $\text{deg}(F)-1$. \\ \\
When $\chi$ is even
\[ (q-1) L(\chi,1)=\left( \sum_{a\in \mathbb{F}_q^*} \chi(a)\right)\sum_{n=0}^{\text{deg}(F)-1} \left(\sum_{D \in \mathcal{M}_n} \chi(D)\right) = \sum_{\substack{D\in \mathbb{F}_q[t] \\ \text{deg}(D)\leq \text{deg}(F)-1}} \chi(D) = 0\]
so a zero is forced at $u=1$. \\ \\
When $\chi$ is odd, we refer to Tao's blog \cite{tao:2019} Theorem 2 where it is stated that all zeros of $L(\chi,u)$ have norm $q^{-1/2}$. \\ \\
We remark that Theorem 9.16B \cite{rosen:2013} only holds for geometric extensions. To each character $\chi$ is associated a cyclic field extension $K_\chi / \mathbb{F}_q[t]$ of degree equals to the order of $\chi$. It is called a geometric extension if its field of constants is $\mathbb{F}_q$. A character associated to a field of constants extension has the form
\[ \chi_c(D) = \zeta^{\text{deg}(D)}\]
where $\zeta$ is a root of unity in $\mathbb{C}$. Those characters aren't Dirichlet characters since they are not periodic to any modulo ($\chi_c$ is never 0), which shows that $K_\chi$ is geometric in the case of Dirichlet characters.
\end{proof}
The following lemma can also be obtained using Poisson summation formula (see \cite{Bui&Florea:2016}). In this case, the character $\chi$ doesn't have to be primitive, but the summands become Gauss sums.
\begin{lemma}[Duality]\label{lm:duality}
Let $\chi$ be a primitive Dirichlet character of conductor $F\neq1$. Then
  \[ \sum_{B\in\mathcal{M}_j} \chi(B) = \omega(\chi) q^{j-\emph{deg}(F)/2} \sum_{k=0}^{\emph{deg}(F)-1-j} \sigma_\chi(k) \sum_{B \in \mathcal{M}_{\emph{deg}(F)-1-j-k}} \overline{\chi}(B) \]
  where $\omega(\chi)$ is given by (\ref{dual:1}) and $\sigma_\chi(k)$ is given by (\ref{dual:2}).
Furthermore, $\sigma_\chi(k)$ only depends on the parity of the character, $|\sigma_\chi(k)| \leq q$, and $|\omega(\chi)|=1$.
\end{lemma}

\begin{proof}
This was proven by Rudnick (\cite{Rudnick:2008} Proposition 7) for quadratic Dirichlet characters. We adjust his computations for arbitrary orders. \\ \\
Assume first that $\chi$ is odd. By the above proposition
\[  L(\chi,u) = \prod_{j=1}^M (1-q^{1/2}e^{i\theta_j}u) \]
where $M = \text{deg}(F)-1$. Since
\[ L(\overline{\chi},u) = \prod_{j=1}^M (1-q^{1/2}e^{-i\theta_j}u) \]
we have
\begin{align}\label{eq:fe}
 L\left(\overline{\chi},\frac{1}{qu}\right) &=  \prod_{j=1}^M \left(1-e^{-i\theta_j} \frac{1}{q^{1/2}u}\right) \\
 &= \frac{\prod_{j=1}^M e^{-i\theta_j}}{(q^{1/2}u)^M} \prod_{j=1}^M (e^{i\theta_j}u q^{1/2}-1)  \nonumber\\
 &= \frac{1}{\omega(\chi) u^M q^{M/2}} L(\chi,u)\nonumber
 \end{align}
 where
 \begin{equation}\label{dual:1} \omega(\chi) = \prod_{j=1}^M -e^{i\theta_j}.
 \end{equation}
  Now, using
 \[ L(\chi,u) = \sum_{j=0}^M \left(\sum_{B\in\mathcal{M}_j} \chi(B)\right) u^j\]
 and the last equation, we get by comparing powers of $u$
 \[ \sum_{B\in\mathcal{M}_j} \chi(B) = \omega(\chi) q^{j-M/2}\sum_{B\in\mathcal{M}_{M-j}} \overline{\chi}(B).\]
 When $\chi$ is even, we have to deal with the extra factor $(1-u)$. Let
 \[ \frac{ L(\chi,u)}{1-u}=\sum_{j=0}^M A_j u^j.\]
Since we removed the extra zero, we can use the above computation to get
 \[  A_j = \omega(\chi) q^{j-M/2}\overline{A}_{M-j}.\]
Also, we have
\[ \sum_{B\in\mathcal{M}_j} \chi(B) =  A_j-A_{j-1}\]
 which implies
 \begin{equation}
 \label{eq:ehx}A_j = \sum_{k=0}^j   \sum_{B\in\mathcal{M}_k} \chi(B).
 \end{equation}
 Then
 \[  \sum_{B\in\mathcal{M}_j} \chi(B) =  A_j-A_{j-1} = \omega(\chi) q^{j-1-M/2}(q\overline{A}_{M-j} - \overline{A}_{M+1-j})\]
 and by Equation \eqref{eq:ehx} above
 \[ \sum_{B\in\mathcal{M}_j} \chi(B) = -\omega(\chi) q^{j-M/2-1} \left( \sum_{B\in\mathcal{M}_{M+1-j}} \overline{\chi}(B) - (q-1) \sum_{k=0}^{M-j} \sum_{B\in\mathcal{M}_k} \overline{\chi}(B) \right) \]
 where $M=\text{deg}(F)-2$. For ease of use we define
\begin{equation}\label{dual:2} \sigma_\chi(k) := \begin{cases}
q^{1/2} & \text{if } k=0\\
 0 & \text{otherwise}
  \end{cases}  \end{equation}
  when $\chi$ is odd and
  \[ \sigma_\chi(k) := \begin{cases}
-1 & \text{if } k=0\\
 q-1 & \text{otherwise}
  \end{cases}  \]
  when $\chi$ is even.
\end{proof}

\section{The explicit formula} \label{sec:comp}

\begin{proposition}[Riemann Hypothesis]\label{prop:rh2}
For $\chi$ a primitive Dirichlet character of conductor $F\neq1$ coprime to $N_E$
\begin{equation}\label{eq:222} L(E\otimes \chi,u)=\prod_{j=1}^{M}  \left(1-qe^{i\theta_{j}}u\right)\end{equation}
where $M=2 \, \emph{deg}(F) + \emph{deg}(N_E)-4$.
\end{proposition}
\begin{proof}
 Let $\ell$ be the order of $\chi$. We first assume $\ell$ to be prime. Let $K_\chi/\mathbb{F}_q[t]$ be the cyclic field extension associated to $\chi$. By Artin's conjecture, which has been proven over function fields, we have
\[ L(E/K_\chi,u) = L(E,u)\prod_{i=1}^{l-1}L(E\otimes\chi^i,u)\]
and all of these functions are entire. When no field is specified in the parameters of the $L$-function, the Euler product is assumed to be over $\mathbb{F}_q[t]$, otherwise it is over the specified field. By [CFJ15] Theorem 1.1 (ii)
\[ L(E/K_\chi,u) =\prod_{j=1}^L \left(1-qe^{i\theta_j}u\right)\]
where $L=\text{deg}_{K_\chi}(\mathcal{N}_E)+2(2g-2)$ and $g$ is the genus of $K_\chi$. The function $\text{deg}_{K_\chi}$ is defined on primes of $K_\chi$ by $\text{deg}_{K_\chi}(\mathfrak{p}):=\text{log}_q |K_\chi/(\mathfrak{p})|$ and extended using $\text{deg}_{K_\chi}(\mathfrak{p}\mathfrak{q})=\text{deg}_{K_\chi}(\mathfrak{p})+\text{deg}_{K_\chi}(\mathfrak{q})$. In particular, for $P$ a prime of $\mathbb{F}_q[t]$, we have $\text{deg}_{K_\chi}(P)=\ell\cdot \text{deg}(P)$. The conductor $\mathcal{N}_E$ is the conductor of $E$ over $K_\chi$. Since $(F,N_E)=1$, only the prime at infinity might ramify in $K_\chi/\mathbb{F}_q[t]$. We recall that we assume $P_\infty^2 \mid N_E$ where $P_\infty$ is the prime at infinity of $\mathbb{F}_q[t]$. \\ \\
If $\chi$ is even, the prime at infinity is not ramified in $K_\chi/\mathbb{F}_q[t]$. This means $\text{deg}_{K_\chi}(\mathcal{N}_E)=\ell\cdot \text{deg}(N_E)$. We also have $2g=(\ell-1)(\text{deg}(F)-2)$ and $L(E,u)$ has degree $\text{deg}(N_E)-4$ (\cite{Brumer:1992}, Appendix). This implies  $M=2 \, \text{deg}(F) + \text{deg}(N_E)-4$ since the twisted $L$-functions have the same number of zeros since $\ell$ is prime. \\ \\
If $\chi$ is odd, the prime at infinity is ramified in $K_\chi/\mathbb{F}_q[t]$. We denote by $\mathfrak{p}_\infty$ the prime at infinity of $K_\chi/\mathbb{F}_q[t]$. We have $P_\infty=(\mathfrak{p}_\infty)^\ell$. However, $\mathfrak{p}_\infty^{2\ell}\nmid \mathcal{N}_E$. We have $\mathfrak{p}_\infty^{2}\mid \mathcal{N}_E$ instead, so the conductors $N_E$ and $\mathcal{N}_E$ are not equal and we have $\text{deg}_{K_\chi}(\mathcal{N}_E)=\ell\cdot \text{deg}(N_E) - 2(\ell -1)$. Also, $2g=(\ell-1)(\text{deg}(F)-1)$.  This implies  $M=2 \, \text{deg}(F) + \text{deg}(N_E)-4$. The result can be generalized to composite orders using induction. \\ \\
We remark that the theorems we have used in this proof include the prime at infinity in their definition of the $L$-function. Since we have assumed that $E$ has additive reduction at infinity, it makes no difference that we exclude the prime at infinity in our definition.
\end{proof}

\begin{lemma}[Explicit Formula]\label{lm:ef}
For $\chi$ a primitive Dirichlet character of conductor $F\neq1$ coprime to $N_E$
\[ \sum_{d\mid n} \sum_{P\in\mathcal{P}_{n/d}}(n/d) (\alpha_P^d + \overline{\alpha}_P^d) \chi^d(P) = -q^n \sum_{j=1}^{M} e^{in\theta_{j}}\]
where $\alpha_P^d+\overline{\alpha}_P^d$ is replaced by $a_P^d$ for primes dividing $N_E$.
\end{lemma}
\begin{proof}
We get the result by comparing the coefficients of the powers of $u$ in the logarithmic derivative of \eqref{eq:222} to those coming from the logarithmic derivative of the Euler product.
\end{proof}

\section{Average of traces}
The main theorem of this paper is the following.
\begin{theorem}\label{thm:main}
For $\epsilon>0$, $n>0$, and $N>4\,\emph{deg}(N_E)$
\begin{align*}
\langle\emph{tr } \Theta^n\rangle_{N,C} &=
\begin{cases}
 1, & \emph{if } n \emph{ is even} \\
 0, & \emph{if } n \emph{ is odd}
  \end{cases} \\
  &+\mathcal{O}_{E,q}\left((n+N)N^{2\,\emph{deg}(N_E)+3}\left( \frac{1}{q^{N/8}} + \frac{1}{q^{\epsilon N}} + \frac{q^{n/2}}{q^{(1-\epsilon)N}}\right) +  n^2q^{-n/4} \right).
\end{align*}
\end{theorem}
\begin{proof}
 By the explicit formula (Lemma \ref{lm:ef})
\begin{equation}\label{eq:main}
\langle\text{tr } \Theta^n\rangle_{N,C} =  \frac{-1}{q^n|\mathcal{H}_{N,C}|} \sum_{d\mid n}\sum_{\text{deg}(P)=n/d} (n/d)(\alpha_P^d+\overline{\alpha}_P^d) \sum_{D \in \mathcal{H}_{N,C}} \chi_D^d(P)
\end{equation}
Combining each estimate of Sections \ref{sc:prime}, \ref{sc:sq}, and \ref{sc:hp} concludes the proof.
\end{proof}

\section{The sieve}
We sieve for two conditions, the polynomials must be square-free and congruent to $C$ modulo $N_E$ for some $C$ coprime to $N_E$. We use the techniques of \cite{Bui&Florea:2016} Lemma 2.2. The notation $A\mid B^\infty$ means that if a prime divides $A$, then it must divide $B$.

\begin{lemma}\label{pr:sieve}
 \[ \sum_{D \in \mathcal{H}_{N,C}}\chi_D(P)= \]
 \[\frac{1}{|(\mathbb{F}_q[t]/(N_E))^*|} \sum_{\psi \emph{ mod } N_E} \overline{\psi}(C)\sum_{k=0}^N\sum_{m=0}^k \alpha_\psi(m)\sum_{\substack{Q_1 \mid M_\psi \\ Q_2 \mid (PN_\psi)^\infty \\ \emph{deg}(Q_1)+2\emph{deg}(Q_2)=k-m}} \mu(Q_1)\hat{\psi}(Q_1)\chi_{Q_1}(P)\hat{\psi_2}(Q_2) \sum_{D\in\mathcal{M}_{N-k}} \hat{\psi}(D)\chi_D(P)\]
 where $\psi_2:=\psi^2$, $\psi=\hat{\psi}\psi_0$ and $\psi_2=\hat{\psi_2}(\psi_2)_0$ as in Equation \eqref{eq:dc}. The modulus of $\psi_0$ is denoted by $M_\psi$ and the modulus of $(\psi_2)_0$ is denoted by $N_\psi$. The function $\alpha_\psi$ is given by \eqref{pr:2} and $|\alpha_\psi(m)|\leq m^{\emph{deg}(N_E)}q^{m/2}$.
\end{lemma}

\begin{proof}
Polynomials $D \in \mathbb{F}_q[t]$ that are congruent to $C$ mod $N_E$ are picked up using
\begin{equation}\label{eq:congr} \frac{1}{|(\mathbb{F}_q[t]/(N_E))^*|} \sum_{\psi \text{ mod } N_E} \overline{\psi}(C)\psi(D)= 
\begin{cases}
 1 & \text{if } D \equiv C \text{ mod } N_E\\
 0 & \text{otherwise}
  \end{cases} 
 \end{equation}
where the sum is over all Dirichlet characters modulo $N_E$. So
\begin{equation}\label{pr:ez}  \sum_{D \in \mathcal{H}_{N,C}}\chi_D(P) =  \frac{1}{|(\mathbb{F}_q[t]/(N_E))^*|} \sum_{\psi \text{ mod } N_E} \overline{\psi}(C)\sum_{D \in \mathcal{H}_{N}}\psi(D)\chi_D(P)\end{equation}
and we write the generating series as
\begin{align}\label{pr:main}
\sum_{h=0}^\infty  \left(\sum_{D \in \mathcal{H}_h}\psi(D)\chi_D(P)\right)u^h &= \prod_{Q \text{ prime} }(1+\psi(Q)\chi_Q(P)u^{\text{deg}(Q)}) \\
&=\frac{\prod_{Q}(1-\psi_2(Q)\chi_Q^2(P) u^{2\text{deg}(Q)})}{\prod_{Q}(1-\psi(Q)\chi_Q(P)u^{\text{deg}(Q)})}\nonumber \\
&=\frac{\prod_{Q}(1-\hat{\psi}(Q)\chi_Q(P)u^{\text{deg}(Q)})^{-1}\prod_{Q\mid M_\psi}(1-\hat{\psi}(Q)\chi_Q(P)u^{\text{deg}(Q)})}{\prod_{Q}(1-\hat{\psi_2}(Q)u^{2\text{deg}(Q)})^{-1}\prod_{Q\mid PN_\psi}(1-\hat{\psi_2}(Q)u^{2\text{deg}(Q)})}. \nonumber
\end{align}
We now expand all four products. We first have
\begin{equation}\label{pr:1}
\prod_{Q}(1-\hat{\psi}(Q)\chi_Q(P)u^{\text{deg}(Q)})^{-1} = \sum_{h=0}^\infty  \left(\sum_{D \in \mathcal{M}_h}\hat\psi(D)\chi_D(P)\right)u^h.
\end{equation}
Also
\[ \prod_{Q}(1-\hat{\psi_2}(Q)u^{2\text{deg}(Q)})^{-1} = L(\hat{\psi_2},u^2).\]
We assume first that $\hat{\psi_2}$ isn't trivial. We use Proposition \ref{prop:rh} to expand
\begin{align*}
 \frac{1}{L(\hat{\psi_2},u^2)} &= (1-u^2)^{-\lambda}\prod_{j=1}^M(1-q^{1/2}e^{i\theta_j,\psi}u^2)^{-1} \\
 &=\left(\sum_{h=0}^\infty u^{2h} \right)^\lambda \prod_{j=1}^M \sum_{h=0}^\infty q^{h/2} e^{ih\theta_j,\psi}u^{2h}.
 \end{align*}
We expand again to get
\begin{equation}\label{pr:2} \frac{1}{L(\hat{\psi_2},u^2)} = \sum_{h=0}^\infty \alpha_\psi(h) u^h\end{equation}
and bounding trivially gives $|\alpha_\psi(h)| \leq h^{\text{deg}(N_E)}q^{h/4}$. \\ \\
If $\hat{\psi_2}$ is trivial, then $1/L(\hat{\psi_2},u^2)=1-qu^2$ and so $\alpha_\psi(0)=1$, $\alpha_\psi(2)=-q$, and it is 0 everywhere else, which is the case of \cite{Bui&Florea:2016}  Lemma 2.2. We use the bound  $|\alpha_\psi(h)| \leq h^{\text{deg}(N_E)}q^{h/2}$ to deal with both cases at once. \\ \\
For the two remaining products, we use
\[ (1-\hat{\psi_2}(Q)u^{2\text{deg}(Q)})^{-1} = \sum_{h=0}^\infty \hat{\psi_2}(Q^h)u^{2h\,\text{deg}(Q)} \]
to get
\begin{equation}\label{pr:3} \frac{\prod_{Q\mid M_\psi}(1-\hat{\psi}(Q)\chi_Q(P)u^{\text{deg}(Q)})}{\prod_{Q\mid PN_\psi}(1-\hat{\psi_2}(Q)u^{2\text{deg}(Q)})} = \sum_{h=0}^\infty\left(\sum_{\substack{Q_1 \mid M_\psi\\ Q_2 \mid (PN_\psi)^\infty \\ \text{deg}(Q_1)+2\text{deg}(Q_2)=h}} \mu(Q_1)\hat{\psi}(Q_1)\chi_{Q_1}(P)\hat{\psi_2}(Q_2) \right)u^h.\end{equation}
We multiply the series (\ref{pr:1}), (\ref{pr:2}), and (\ref{pr:3}) and we compare powers of $u$ in Equation \eqref{pr:main}. We put the result into Equation \eqref{pr:ez} to conclude.
\end{proof}

We now compute the size of $\mathcal{H}_{N,C}$ in a similar fashion, but we use Perron's formula instead of multiplying the generating series together.
\begin{lemma}\label{lm:size}
For $(C,N_E)=1$ and any $\epsilon>0$
\[   |\mathcal{H}_{N,C}| = \frac{1}{|(\mathbb{F}_q[t]/(N_E))^*|}(1-q^{-1})\prod_{\substack{Q \mid N_E \\ Q \emph{ prime}}}(1+q^{-\emph{deg}(Q)})^{-1} q^{N+1}+\mathcal{O}_{E,q,\epsilon}(q^{N(1/4+\epsilon)}). \]
In particular, the family $\mathcal{H}^*_N$ is equidistributed in the invertible congruence classes modulo $N_E$ as $N\rightarrow\infty$. Furthermore, $|\mathcal{H}_{N,C}| \asymp_{E,q}  q^{N}$.
\end{lemma}
\begin{proof}
We deal with the congruence condition using \eqref{eq:congr}.
 So
\begin{equation}\label{pr:ez2}  \sum_{D \in \mathcal{H}_{N,C}}1 =  \frac{1}{|(\mathbb{F}_q[t]/(N_E))^*|} \sum_{\psi \text{ mod } N_E} \overline{\psi}(C)\sum_{D \in \mathcal{H}_{N}}\psi(D)
\end{equation}
and we write the generating series as
\begin{align*}
\mathcal{L}(\psi,u):=\sum_{h=0}^\infty  \left(\sum_{D \in \mathcal{H}_h}\psi(D)\right)u^h &= \prod_{Q \text{ prime} }(1+\psi(Q)u^{\text{deg}(Q)}) \\
&=\frac{\prod_{Q}(1-\psi_2(Q) u^{2\text{deg}(Q)})}{\prod_{Q}(1-\psi(Q)u^{\text{deg}(Q)})} \\
&=\frac{\prod_{Q}(1-\hat{\psi}(Q)u^{\text{deg}(Q)})^{-1}\prod_{Q\mid M_\psi}(1-\hat{\psi}(Q)u^{\text{deg}(Q)})}{\prod_{Q}(1-\hat{\psi_2}(Q)u^{2\text{deg}(Q)})^{-1}\prod_{Q\mid N_\psi}(1-\hat{\psi_2}(Q)u^{2\text{deg}(Q)})} \\
&=\frac{L(\hat{\psi},u)\prod_{Q\mid M_\psi}(1-\hat{\psi}(Q)u^{\text{deg}(Q)})}{L(\hat{\psi}_2,u^2)\prod_{Q\mid N_\psi}(1-\hat{\psi_2}(Q)u^{2\text{deg}(Q)})}
\end{align*}
where $\psi_2:=\psi^2$, $\psi=\hat{\psi}\psi_0$ and $\psi_2=\hat{\psi_2}(\psi_2)_0$ as in Equation \eqref{eq:dc}. The modulus of $\psi_0$ is denoted by $M_\psi$ and the modulus of $(\psi_2)_0$ is denoted by $N_\psi$. Perron's formula works by dividing the generating series by $u^{N+1}$ in order to create a pole at $u=0$ such that its residue is the $N$th coefficient of the series. Here
\[ \frac{\mathcal{L}(\psi,u)}{n^{N+1}} =  \sum_{h=0}^\infty  \left(\sum_{D \in \mathcal{H}_h}\psi(D)\right)u^{h-N-1}.\]
This sum is the Laurent expansion of $\frac{\mathcal{L}(\psi,u)}{n^{N+1}}$ at $u=0$, so the residue at $u=0$ is the coefficient of $u^{-1}$, which is
\[ \sum_{D \in \mathcal{H}_N}\psi(D).\]
 We're going to use Cauchy's residue theorem by integrating on the complex circle $\mathcal{C}_\epsilon: \,|u|=q^{-1/4-\epsilon}$ the following function
\begin{equation}\label{eq:def}
\frac{ \mathcal{L}(\psi,u)}{u^{N+1}}= \frac{L(\hat{\psi},u)\prod_{Q\mid M_\psi}(1-\hat{\psi}(Q)u^{\text{deg}(Q)})}{L(\hat{\psi}_2,u^2)\prod_{Q\mid N_\psi}(1-\hat{\psi_2}(Q)u^{2\text{deg}(Q)})} \frac{1}{u^{N+1}}.
 \end{equation}
If $\hat{\psi}$ is trivial
\[ L(\hat{\psi},u) = \frac{1}{1-qu}\]
and it has a unique simple pole at $u=1/q$. If $\hat{\psi}$ isn't trivial, then
\[  L(\hat{\psi},u)  = (1-u)^{\lambda_\psi} \prod_{j=1}^{K_\psi} (1-q^{1/2}e^{i\theta_{\psi,j}}u) \]
where $K_\psi \leq \text{deg}(N_E)$ and it has no poles. If $\hat{\psi}_2$ is trivial
\[ \frac{1}{L(\hat{\psi}_2,u^2)} = 1-qu^2\]
and it has no poles. If $\hat{\psi}_2$ isn't trivial
\[ \frac{1}{ L(\hat{\psi}_2,u^2)}  = (1-u)^{-\lambda_{\psi_2}} \prod_{j=1}^{K_{\psi_2}} (1-q^{1/2}e^{i\theta_{\psi_2,j}}u^2)^{-1} \]
where $K_{\psi_2} \leq \text{deg}(N_E)$ and it has poles of norm $q^{-1/4}$. The two remaining products of (\ref{eq:def}) have poles and zeros of norm one. \\ \\
 We define
\[ A(E,q,\epsilon):=\text{sup}_{\psi \text{ mod } N_E} \text{sup}_{|u|=q^{-1/4-\epsilon}}\left| \frac{L(\hat{\psi},u)\prod_{Q\mid M_\psi}(1-\hat{\psi}(Q)u^{\text{deg}(Q)})}{L(\hat{\psi}_2,u^2)\prod_{Q\mid N_\psi}(1-\hat{\psi_2}(Q)u^{2\text{deg}(Q)})}\right|.\]
We have $A(E,q,\epsilon)<\infty$  since these functions have no poles on $\mathcal{C}_\epsilon$. We also have $|1/u^{N+1}| = q^{N(1/4+\epsilon)+1/4+\epsilon}$ on $\mathcal{C}_\epsilon$, so
\begin{equation}\label{eq:tard} \left| \oint_{\mathcal{C}_\epsilon} \frac{\mathcal{L}(\psi,u)}{u^{N+1}} \right| \leq 2\pi A(E,q,\epsilon) q^{N(1/4+\epsilon)} \ll_{E,q,\epsilon} q^{N(1/4+\epsilon)}.\end{equation}
If $\hat{\psi}$ is trivial, then $\psi$ is the principal character modulo $N_E$. In this case, we have by Cauchy's residue theorem
\begin{align*} \frac{1}{2\pi i} \oint_{\mathcal{C}_\epsilon}\frac{\mathcal{L}(\psi,u)}{u^{N+1}} &=  \text{Res}_{u=0} \frac{\mathcal{L}(\psi,u)}{u^{N+1}}  + \text{Res}_{u=1/q} \frac{\mathcal{L}(\psi,u)}{u^{N+1}} \\
&=\sum_{D \in \mathcal{H}_N} \psi(D) + (1-q^{-1})\prod_{\substack{Q \mid N_E \\ Q \text{ prime}}}(1+q^{-\text{deg}(Q)})^{-1} q^{N+1}
\end{align*}
since $\hat{\psi}_2$ is also trivial, and $M_{\psi}$ and $N_\psi$ are both equal to the product of all the primes dividing $N_E$. \\
If $\psi$ isn't the principal character modulo $N_E$, then there is only a pole at $u=0$ and
\[ \frac{1}{2\pi i} \oint_{\mathcal{C}_\epsilon}\frac{\mathcal{L}(\psi,u)}{u^{N+1}} =  \sum_{D \in \mathcal{H}_N} \psi(D).\]
We use those two results along with \eqref{eq:tard} in \eqref{pr:ez2} to get
\[  \sum_{D \in \mathcal{H}_{N,C}}1 =  \frac{1}{|(\mathbb{F}_q[t]/(N_E))^*|}(1-q^{-1}) q^{N+1}\prod_{\substack{Q \mid N_E \\ Q \text{ prime}}}(1+q^{-\text{deg}(Q)})^{-1}+\mathcal{O}_{E,q,\epsilon}(q^{N(1/4+\epsilon)}).\]
By Proposition 1.4. and 1.6. \cite{rosen:2013} we have $|(\mathbb{F}_q[t]/(N_E))^*|\asymp q^{\text{deg}(N_E)}$, so we can conclude
\[ |\mathcal{H}_{N,C}| \asymp_{E,q} q^N.\]
\end{proof}

\section{Contribution of the primes}\label{sc:prime}
The contribution of the primes in Equation \eqref{eq:main} corresponds to the terms with $d=0$, which is
 \[ S_C(N,n):= \frac{-1}{q^n|\mathcal{H}_{N,C}|} \sum_{P\in\mathcal{P}_{n}} n(\alpha_P+\overline{\alpha}_P) \sum_{D \in \mathcal{H}_{N,C}} \chi_D(P). \]
We start with a lemma that will be useful for bounding some quantities coming from the sieve.
 \begin{lemma}\label{lm:ind}
 \[ \sum_{\substack{Q\mid N_E^\infty \\ \emph{deg}(Q)=N}}  1  \leq  (N+1)^{\emph{deg}(N_E)}\]
 \end{lemma}
 \begin{proof}
 In the worst case scenario, $N_E$ is a product of distinct primes of degree one. We use induction on the degree of $N_E$. If $N_E$ is prime, then there is only the possibility $Q=N_E^N$, the lemma is true in this case. Now we assume the lemma is true for $\text{deg}(N_E)=k$. Let $\text{deg}(N_E)=k+1$ and fix a prime $P_0$ that divides $N_E$. We split the terms depending on the powers of $P_0$ dividing them
 \[\sum_{\substack{Q\mid N_E^\infty \\ \text{deg}(Q)=N}} 1 = \sum_{j=0}^N \sum_{\substack{Q\mid (N_E/P_0)^\infty \\ \text{deg}(Q)=N-j}} 1 \leq (N+1)^{k+1} \]
 and this concludes the induction.
 \end{proof}
 
\begin{proposition}
For any $\epsilon>0$ and $N>4\,\emph{deg}(N_E)$
\[  S_C(N,n) \ll_{E,q} (n+N) N^{2\,\emph{deg}(N_E)+3}\left( \frac{1}{q^{N/8}} + \frac{1}{q^{\epsilon N}} + \frac{q^{n/2}}{q^{(1-\epsilon)N}}\right).\]
\end{proposition}

\begin{proof}

We apply the sieve (Proposition \ref{pr:sieve}) and quadratic reciprocity $\chi_D(P)=(-1)^{\text{deg}(D)\text{deg}(P)(q-1)/2}\chi_P(D)$ to write
  \begin{align}\label{eq:cpmain}
  S_C(N,n)=\frac{-1}{q^n|\mathcal{H}_{N,C}|} \sum_{P\in\mathcal{P}_{n}} n(\alpha_P+\overline{\alpha}_P)\frac{1}{|(\mathbb{F}_q[t]/(N_E))^*|} \sum_{\psi \text{ mod } N_E} \overline{\psi}(C)\sum_{k=0}^N\sum_{m=0}^k \alpha_\psi(m) \sum_{\ell=0}^{\lfloor \frac{k-m}{2n} \rfloor}  \\
 \sum_{\substack{Q_1 \mid M_\psi \\ Q_2 \mid \tilde{N}_\psi^\infty \\ \text{deg}(Q_1)+2\text{deg}(Q_2)=k-m-2n\ell}} (-1)^{n(N-k)(q-1)/2}\mu(Q_1)\hat{\psi}(Q_1)\chi_{Q_1}(P)\hat{\psi_2}(Q_2P^\ell) \sum_{D\in\mathcal{M}_{N-k}} \hat{\psi}\chi_P(D) \nonumber
  \end{align}
  where we split the sum over $Q_2$ according to powers of $P$, so we define $\tilde{N}_\psi:=N_\psi/P$ if $P\mid N_\psi$, and $\tilde{N}_\psi:=N_\psi$ if $P\nmid N_\psi$.
  Assume first that $n\leq N/4$. The degree of the conductor of $\hat{\psi}\chi_P$ is then bounded by $N/4+\text{deg}(N_E)$. Then, when $N-k \geq N/4+\text{deg}(N_E)$, the sum over the monic polynomials is zero, so it is bounded by $q^{N/4+\text{deg}(N_E)}$. Using Lemma \ref{lm:size} and Lemma \ref{lm:ind} and noticing $Q_1$ has at most $2^{\text{deg}(N_E)}$ possibilities gives
  \begin{equation}\label{bound:1} S_C(N,n) \ll_{E,q} \frac{N^{2\,\text{deg}(N_E)+3}}{q^{N/8}}\end{equation}
  by bounding everything else trivially.
 Therefore, the contribution of the primes of low degrees tends to zero as $N\rightarrow\infty$. \\ \\
Now, assume $n>N/4$. We also assume $N>4\,\text{deg}(N_E)$ so that $\hat{\psi}\chi_P$ is primitive. We split into two cases depending on whether $k\geq \epsilon N$ or $k<\epsilon N$ for any $\epsilon>0$. \\ \\
We start with the case $k\geq \epsilon N$. We want to bring the sum over $P$ inside in order to use the explicit formula (Lemma \ref{lm:ef}). The sum is then
 \begin{align}\label{eq:lowk}
  \frac{-1}{q^n|\mathcal{H}_{N,C}|} \frac{1}{|(\mathbb{F}_q[t]/(N_E))^*|} \sum_{\psi \text{ mod } N_E} \overline{\psi}(C)\sum_{k\geq\epsilon N}^N\sum_{m=0}^k \alpha_\psi(m)  \sum_{\ell=0}^{\lfloor \frac{k-m}{2n} \rfloor} \sum_{\substack{Q_1 \mid M_\psi \\ Q_2 \mid \tilde{N}_\psi^\infty \\ \text{deg}(Q_1)+2\text{deg}(Q_2)=k-m-2n\ell}} \\
 \mu(Q_1)\hat{\psi}(Q_1)\hat{\psi_2}(Q_2) \sum_{D\in\mathcal{M}_{N-k}}\hat{\psi}(D) \sum_{P\in\mathcal{P}_{n}} n(\alpha_P+\overline{\alpha}_P)\chi_D(P)\chi_{Q_1}(P)\hat{\psi_2}(P^\ell) \nonumber
  \end{align}
where we reapplied quadratic reciprocity on $\chi_P(D)$. Let $\Psi:=\chi_{DQ_1}\hat{\psi}^\ell_2$ and $\Psi=\hat{\Psi}\Psi_0$ as in (\ref{eq:dc}). The explicit formula for $L(E\otimes\hat{\Psi},u)$ gives the bound
\begin{equation*}
\sum_{d\mid n} \sum_{P\in\mathcal{P}_{n/d}} (n/d)(\alpha_P^d+\overline{\alpha}_P^d)\hat{\Psi}^d(P) \ll_E N q^n.
 \end{equation*}
It is important to mention here that the conductor of $\hat{\Psi}$ might not be coprime to $N_E$. The Riemann hypothesis is still valid for $L(E\otimes\hat{\Psi},u)$, but the number of zeros given in Proposition \ref{prop:rh2} could be changed by a quantity bounded by $2\,\text{deg}(N_E)$, which is why the bound remains valid.  This bound also holds when $\hat{\Psi}$ is trivial, since then we get $L(E,u)$ which has $ \text{deg}(N_E)-4$ zeros of norm $1/q$.
Now, since
 \begin{align*}
  \sum_{d\mid n} \sum_{P\in\mathcal{P}_{n/d}} (n/d)(\alpha_P^d+\overline{\alpha}_P^d)\hat{\Psi}^d(P)  &=  \\
  \sum_{P\in\mathcal{P}_{n}} n(\alpha_P+\overline{\alpha}_P)\hat{\Psi}(P) &+ \sum_{\substack{d\mid n \\ d>1}} \sum_{P\in\mathcal{P}_{n/d}} (n/d)(\alpha_P^d+\overline{\alpha}_P^d)\hat{\Psi}^d(P)
  \end{align*}
  and
  \[ \sum_{\substack{d\mid n \\ d>1}} \sum_{P\in\mathcal{P}_{n/d}} (n/d)(\alpha_P^d+\overline{\alpha}_P^d)\hat{\Psi}^d(P) = \mathcal{O}(n q^n)\]
   we have
  \begin{equation}\label{eq:weil}
    \sum_{P\in\mathcal{P}_{n}} n(\alpha_P+\overline{\alpha}_P)\hat{\Psi}(P)  \ll_E q^n(N+n). 
    \end{equation}
Also, we have
\[  \sum_{P\in\mathcal{P}_{n}} n(\alpha_P+\overline{\alpha}_P)\Psi(P) =  \sum_{P\in\mathcal{P}_{n}} n(\alpha_P+\overline{\alpha}_P)\hat{\Psi}(P) -  \sum_{\substack{P\in\mathcal{P}_{n} \\ P \mid K_\Psi}} n(\alpha_P+\overline{\alpha}_P)\hat{\Psi}(P) \]
where $K_\Psi$ is the modulus of $\Psi_0$. We have assumed $n>N/4$ and $N>4\,\text{deg}(N_E)$, so only $\chi_D$ can contribute to $K_\Psi$, and it contributes at most two different primes. So
\[  \sum_{P\in\mathcal{P}_{n}} n(\alpha_P+\overline{\alpha}_P)\Psi(P) =  \sum_{P\in\mathcal{P}_{n}} n(\alpha_P+\overline{\alpha}_P)\hat{\Psi}(P) + \mathcal{O}(nq^{n/2}). \]
Then
\[  \sum_{P\in\mathcal{P}_{n}} n(\alpha_P+\overline{\alpha}_P)\Psi(P) \ll_E q^n(N+n).\]
This sum is exactly the sum over $\mathcal{P}_n$ in (\ref{eq:lowk}). We use $|\alpha_\psi(m)|\leq q^{k/2}N^{\text{deg}(N_E)}$ from Lemma \ref{pr:sieve} since $m\leq k$ and we bound everything else trivially as we did for (\ref{bound:1}). We get that (\ref{eq:lowk}) is bounded by
\begin{equation}\label{bound:2} \ll_{E,q} \frac{(n+N)N^{2\text{deg}(N_E)+3}}{q^{\epsilon N/2}}\end{equation}
which goes to zero as $N\rightarrow \infty$. \\ \\
For the case $k<\epsilon N$, we apply duality (Lemma \ref{lm:duality}) to \eqref{eq:cpmain} to get
    \begin{align}\label{eq:eh}
  \frac{-1}{q^n|\mathcal{H}_{N,C}|} \sum_{P\in\mathcal{P}_{n}} n(\alpha_P+\overline{\alpha}_P)\frac{1}{|(\mathbb{F}_q[t]/(N_E))^*|} \sum_{\psi \text{ mod } N_E} \overline{\psi}(C)\sum_{k=0}^{\lfloor \epsilon N\rfloor}\sum_{m=0}^k \alpha_\psi(m) \sum_{\ell=0}^{\lfloor \frac{k-m}{2n} \rfloor}\\
 \sum_{\substack{Q_1 \mid M_\psi \\ Q_2 \mid \tilde{N}_\psi^\infty \\ \text{deg}(Q_1)+2\text{deg}(Q_2)=k-m-2n\ell}}
 (-1)^{n(N-k)(q-1)/2} \mu(Q_1) \hat{\psi}(Q_1)\chi_{Q_1}(P)\hat{\psi_2}(Q_2P^\ell) \nonumber \\
\omega(\hat{\psi}\chi_P) q^{N-k-(n+\text{deg}(C_{\hat{\psi}}))/2}
 \sum_{r=0}^{n+\text{deg}(C_{\hat{\psi}})-1-N+k} \sigma_{\hat{\psi}\chi_P}(r) \sum_{D \in \mathcal{M}_{n+\text{deg}(C_{\hat{\psi}})-1-N+k-r}} \overline{\hat{\psi}\chi_P}(D) \nonumber
  \end{align}
  where $C_{\hat{\psi}}$ is the conductor of $\hat{\psi}$. \\ \\
The goal is again to sum over $\mathcal{P}_n$ first in order to use the explicit formula (Lemma \ref{lm:ef}). To bring the sum inside, we must deal with $\omega(\hat{\psi}\chi_P)$ to remove its dependency on $P$.  \\ \\
Definitions related to Gauss sums can be found at the end of Section \ref{sec:background}. Corollary 2.4 \cite{David&Florea&Lalin:2019} states that for primitive characters of conductor $F$
  \[ \omega(\chi)= 
\begin{cases}
\frac{1}{\tau(\chi)} q^{-(\text{deg}(F)-1)/2}G(\chi) & \text{if } \chi \text{ odd} \\
q^{-\text{deg}(F)/2}G(\chi) & \text{if } \chi \text{ even}
  \end{cases} 
 \]
 and adjusting Lemma 2.12 (i) \cite{David&Florea&Lalin:2019} for general Dirichlet characters gives
\[ G(\hat{\psi}\chi_P) = \hat{\psi}(P)\chi_P(C_{\hat{\psi}})G(\hat{\psi})G(\chi_P)\]
because $\text{deg}(P)$ is large enough for $P$ to be coprime to $C_{\hat{\psi}}$. Then, assuming the most complicated case where all characters are odd
\begin{align*} \omega(\hat{\psi}\chi_P)&=\frac{1}{\tau(\hat{\psi}\chi_P)}q^{-(\text{deg}(C_{\hat{\psi}})+n-1)/2} G(\hat{\psi}\chi_P) \\
&=\frac{1}{\tau(\hat{\psi}\chi_P)}q^{-(\text{deg}(C_{\hat{\psi}})+n-1)/2}\hat{\psi}(P)\chi_P(C_{\hat{\psi}})G(\hat{\psi})G(\chi_P) \\
&=\frac{1}{\tau(\hat{\psi}\chi_P)}q^{-1/2}\hat{\psi}(P)\chi_P(C_{\hat{\psi}})\omega(\hat{\psi})\tau(\hat{\psi})\omega(\chi_P)\tau(\chi_P).
 \end{align*}
If $\text{deg}(P)$ is even, then $\chi_P(a)=1$ for all $a \in \mathbb{F}_q^*$ and if $\text{deg}(P)$ is odd, $\chi_P(a)=1$ if $a$ is a square in $\mathbb{F}_q^*$ and $\chi_P(a)=-1$ otherwise. This implies $\tau(\hat{\psi}\chi_P)$ and $\tau(\chi_P)$ does not depend on the actual $P$, only on its degree, since the sum defining these quantities is over $\mathbb{F}_q^*$. We recall $\omega(\chi_P)=1$ since $\chi_P$ is quadratic. This implies
\[ \omega_n(\hat{\psi}):=\frac{\omega(\hat{\psi}\chi_P)}{\hat{\psi}(P)\chi_P(C_{\hat{\psi}})}\]
only depends on $\text{deg}(P)=n$, hence the notation. \\ \\
The other quantity that might depend on $P$ is $\sigma_{\hat{\psi}\chi_P}(r)$. We recall that it only depends on the parity of the character, and by the discussion above, the parity of $\hat{\psi}\chi_P$ only depends on $\hat{\psi}$ and the degree of $P$. We therefore use the notation $\sigma_{\hat{\psi},n}(r):=\sigma_{\hat{\psi}\chi_P}(r)$.\\ \\
Replacing in \eqref{eq:eh} gives
    \begin{align}\label{eq:last}
  \frac{-1}{q^n|\mathcal{H}_{N,C}|}  \frac{1}{|(\mathbb{F}_q[t]/(N_E))^*|} \sum_{\psi \text{ mod } N_E} \overline{\psi}(C)\sum_{k=0}^{\lfloor \epsilon N\rfloor}\sum_{m=0}^k \alpha_\psi(m)\sum_{\ell=0}^{\lfloor\frac{k-m}{2n} \rfloor}\sum_{\substack{Q_1 \mid M_\psi \\ Q_2 \mid \tilde{N}_\psi^\infty \\ \text{deg}(Q_1)+2\text{deg}(Q_2)=k-m-2n\ell}} \\
\mu(Q_1)\hat{\psi}(Q_1)\hat{\psi_2}(Q_2)\omega_n(\hat{\psi}) q^{N-k-(n+\text{deg}(C_{\hat{\psi}}))/2}\sum_{r=0}^{n+\text{deg}(C_{\hat{\psi}})-1-N+k} (-1)^{n(n-r-1)(q-1)/2}\nonumber\\ \sigma_{\hat{\psi},n}(r) \sum_{D \in \mathcal{M}_{n+\text{deg}(C_{\hat{\psi}})-1-N+k-r}} \overline{\hat{\psi}}(D) \sum_{P\in\mathcal{P}_{n}} n(\alpha_P+\overline{\alpha}_P) \chi_D(P)\hat{\psi}(P)\chi_{C_{\hat{\psi}}}(P)\chi_{Q_1}(P)\hat{\psi_2}(P^\ell)\nonumber
  \end{align}
after applying quadratic reciprocity on $\chi_P(D)$ and $\chi_P(C_{\hat{\psi}})$. We use (\ref{eq:weil}) to bound the sum over $\mathcal{P}_n$. We have an extra $\hat{\psi}\chi_{C_{\hat{\psi}}}$ in the character, so it adds at most $2\,\text{deg}(N_E)$ to the degree of the conductor. Bounding everything else trivially, we have that (\ref{eq:last}) is
\[ \ll_{E,q} (n+N)N^{2\text{deg}(N_E)+3} q^{n/2}q^{-N+\epsilon N/2}.\]
Combining this bound with (\ref{bound:1}) and (\ref{bound:2}) concludes to proof.
\end{proof}

\section{Contribution of the squares}\label{sc:sq}
The contribution of the squares comes from the terms with even $d$ in Equation \eqref{eq:main}.
\begin{proposition}
\begin{equation*}
\frac{-1}{q^n|\mathcal{H}_{N,C}|}\sum_{\substack{d\mid n \\ 2 \mid d}}\sum_{\emph{deg}(P)=n/d} (n/d)(\alpha_P^d+\overline{\alpha}_P^d) \sum_{D \in \mathcal{H}_{N,C}} \chi_D^d(P)=
\begin{cases}
1 + \mathcal{O}_E\left( n^2q^{-n/4}\right),  &\text{if } 2\,|\,n \\
0,   &\text{if } 2\nmid n
\end{cases} 
\end{equation*}
where $\alpha_P^d+\overline{\alpha}_P^d$ is replaced by $a_P^d$ for bad primes.
\end{proposition}
\begin{proof}
If $2 \nmid n$, there are no terms, hence no contribution. \\ \\ When $2\mid n$, we rewrite the sum as
\begin{equation}\label{eq:squares}
 \frac{-1}{q^{2m}|\mathcal{H}_{N,C}|}\sum_{d|m}\sum_{\text{deg}(P)=m/d} (m/d)\left(\alpha_P^{2d} + \overline{\alpha_P}^{2d}\right)\sum_{D \in \mathcal{H}_{N,C}} \chi_D^{2d}(P)
\end{equation}
where $m:=n/2$. The character equals $\mathds{1}_{P\nmid D}$ because its power is even. A simple sieving gives
\[
\sum_{\substack{D \in \mathcal{H}_{N,C} \\ (D,P)=1}} 1 = \sum_{j=0}^{\lfloor N/\text{deg}(P)\rfloor} (-1)^j \sum_{D \in \mathcal{H}_{N-j\,\text{deg}(P),CP^{-j}}} 1
\]
if $(P,N_E)=1$ so that $P$ is invertible modulo $N_E$. If $P \mid N_E$, then
\[ \sum_{\substack{D \in \mathcal{H}_{N,C} \\ (D,P)=1}} 1  = |\mathcal{H}_{N,C}|\]
since every $D$ is coprime to $N_E$ because $(C,N_E)=1$.
In the case $(P,N_E)=1$, the sum with $j=0$ is exactly $|\mathcal{H}_{N,C}|$. When $j>0$
\[  \sum_{j=1}^{\lfloor N/\text{deg}(P)\rfloor} (-1)^j \sum_{D \in \mathcal{H}_{N-j\,\text{deg}(P),CP^{-j}}} 1 \asymp  \sum_{j=1}^{\lfloor N/\text{deg}(P)\rfloor} (-1)^j q^{N-j\text{deg}(P)} \ll q^{N-\text{deg}(P)} \sum_{j=0}^\infty q^{-j} \ll q^{N-\text{deg}(P)}.\]
Then
\begin{equation}\label{eq:charsum}
\sum_{D \in \mathcal{H}_{N,C}} \chi_D^{2d}(P) = |\mathcal{H}_{N,C}| + \mathcal{O}(q^{N-\text{deg}(P)})
\end{equation}
which holds in both cases. \\  \\
Now (\ref{eq:squares}) is equal to
\begin{equation*}
 \frac{-1}{q^{2m}}\sum_{d|m}\sum_{\text{deg}(P)=m/d}(m/d)\left(\alpha_P^{2d} + \overline{\alpha_P}^{2d}\right)+ \mathcal{O} \left(m^2 q^{-m}\right).
\end{equation*}
We retrieve the double sum by looking at the symmetric square of $L(E,u)$, which is defined by
\begin{align*}
 L(\text{Sym}^2 E, u):=&\prod_{P \nmid N_E} \left( 1-\alpha_P^2 u^{\text{deg}(P)}\right)^{-1} \left( 1-\alpha_P\overline{\alpha}_P u^{\text{deg}(P)}\right)^{-1} \left( 1-\overline{\alpha}_P^2 u^{\text{deg}(P)}\right)^{-1}\\
&\prod_{P\mid N_E}(1-a_P^2 u^{\text{deg}(P)})^{-1} .
\end{align*}
It is related to the variety
\[ \text{Sym}^2 E: y^2 = (x^3+Ax+B)(z^3+Az+B).\]
It is known that (\cite{Cha&Fiorilli&Jouve:2015} Theorem 1.1)
\[ L(\text{Sym}^2 E,u) = \prod_{j=1}^M (1-q^{3/2}e^{i\theta_{\text{Sym}^2 E,j}} u) \]
for some $M<\infty$ depending on $E$. The explicit formula for $L(\text{Sym}^2 E,u)$ is then
\[  \sum_{d|m}\sum_{\text{deg}(P)=m/d}(m/d)\left(\alpha_P^{2d} + q^{d \cdot \text{deg}(P)} + \overline{\alpha_P}^{2d}\right) = q^{3m/2} \sum_{j=1}^M e^{im\theta_{\text{Sym}^2 E,j}}\]
where we replace $\alpha_P^{2d} + q^{d \cdot \text{deg}(P)} + \overline{\alpha_P}^{2d}$ by $a_P^{2d}$ for bad primes. We have
\[ \sum_{d|m}\sum_{\text{deg}(P)=m/d}(m/d)  q^{d \cdot \text{deg}(P)} = q^{2m} \]
and
\[ \sum_{d|m}\sum_{\text{deg}(P)=m/d}(m/d)  q^{d \cdot \text{deg}(P)} = \sum_{d|m}\sum_{\substack{\text{deg}(P)=m/d\\ P \text{ good}}}(m/d)  q^{d \cdot \text{deg}(P)} +\sum_{d|m}\sum_{\substack{\text{deg}(P)=m/d \\ P \text{ bad}}}(m/d)  q^{d \cdot \text{deg}(P)}.\]
Since there are at most $\text{deg}(N_E)$ bad primes
\[ \sum_{d|m}\sum_{\substack{\text{deg}(P)=m/d \\ P \text{ bad}}}(m/d)  q^{d \cdot \text{deg}(P)} \ll \text{deg}(N_E) m q^m.\]
This implies
\[  \sum_{d|m}\sum_{\text{deg}(P)=m/d} (m/d)\left(\alpha_P^{2d} + \overline{\alpha_P}^{2d}\right) = -q^{2m} + \mathcal{O}_E (q^{3m/2})\]
where we replace $\alpha_P^{2d}+ \overline{\alpha_P}^{2d}$ by $a_P^{2d}$ for bad primes and this concludes the proof.
\end{proof}

\section{Contribution of higher powers}\label{sc:hp}
\begin{proposition}
\begin{equation*}
\frac{1}{q^n|\mathcal{H}_{N,C}|}\sum_{\substack{d\mid n \\ d>2}}\sum_{\emph{deg}(P)=n/d} (n/d)(\alpha_P^d+\overline{\alpha}_P^d) \sum_{D \in \mathcal{H}_{N,C}} \chi_D^d(P)\ll n^2 q^{-n/6}.
\end{equation*}
\end{proposition}

\begin{proof}
We bound everything trivially.
\end{proof}

\section{One-level density}\label{sc:old}
\begin{corollary}\label{cor:1ld}
For $\phi \in \mathcal{S}(\mathbb{R})$ an even function such that \emph{supp}$(\hat{\phi})\subset (-1,1)$, we have
\[ \langle Z_{\phi}\rangle_{N,C}=\int_{O(M)} Z_{\phi}(\Theta) d\Theta + \mathcal{O}_{E,q}(1/N)\]
where $M=2N+\emph{deg}(N_E)-4$.
\end{corollary}
\begin{proof}
 The Fourier expansion of $Z_{\phi}(\Theta)$ is
\[ Z_{\phi}(\Theta) = \int_{\mathbb{R}} \phi(x)dx + \frac{1}{M} \sum_{n\neq 0}\hat{\phi}\left(\frac{n}{M}\right) \text{tr } \Theta^n.\]
Since $\text{tr }\Theta^{-n}= \overline{\text{tr }\Theta^n}$, our main result is valid when $n<0$.
Averaging over $\mathcal{H}_{N,C}$ and applying Theorem \ref{thm:main} gives
\begin{align*} \langle Z_{\phi}\rangle_{N,C}&=\hat{\phi}(0)+\frac{1}{M}\sum_{n\neq 0} \hat{\phi}\left(\frac{n}{M}\right) \eta_n \\
&+ \mathcal{O}_{E,q}\left(\frac{1}{M}\sum_{n=1}^\infty \hat{\phi}\left(\frac{n}{M}\right)\left((n+N)N^{2\,\text{deg}(N_E)+3}\left( \frac{1}{q^{N/8}} + \frac{1}{q^{\epsilon N}} + \frac{q^{n/2}}{q^{(1-\epsilon)N}}\right) +  n^2q^{-n/4}\right)\right)
\end{align*}
where
\begin{align*}
\eta_n&=
\begin{cases}
 1, & \text{if } n \text{ is even} \\
 0, & \text{if } n \text{ is odd.}
  \end{cases}
\end{align*}
By Equation \eqref{eq:ortho}
\begin{equation}\label{eq:rank}
 \hat{\phi}(0)+\frac{1}{M}\sum_{n\neq 0} \hat{\phi}\left(\frac{n}{M}\right) \eta_n  = \int_{O(M)} Z_{\phi}(\Theta) d\Theta.
 \end{equation}
For the error term, we have that
\[ (n+N)N^{2\,\text{deg}(N_E)+3}\left( \frac{1}{q^{N/8}} + \frac{1}{q^{\epsilon N}} + \frac{q^{n/2}}{q^{(1-\epsilon)N}}\right)  \]
tends to zero as $N\rightarrow \infty$ provided that $n<2(1-2\epsilon)N$. The range of $n$ is controlled by restricting the support of $\hat{\phi}$, so it must be inside $(-1,1)$ since $M\sim 2N$ as $N\rightarrow\infty$. For the second term of the error, notice that
\[ \sum_{n=1}^\infty n^2q^{-n/4}\]
converges, so
\[\frac{1}{M}  \sum_{n=1}^\infty \hat{\phi}\left(\frac{n}{M}\right) n^2q^{-n/4} \]
tends to zero as $N\rightarrow \infty$.
\end{proof}
\textbf{Remark.} By (6) \cite{Young:2004}, we cannot distinguish between the symmetry types $O$, $SO(\text{even})$, and $SO(\text{odd})$ since they all have the same one-level density when the support of $\hat{\phi}$ is contained inside $(-1,1)$.

\section{Average rank and non-vanishing}\label{sc:av}
One application of the one-level density is getting a bound on the average analytic rank. 
\begin{theorem}\label{thm:rank}
The family $\mathcal{H}_{N,C}$ of quadratic twists of an elliptic curve over $\mathbb{F}_q[t]$ has average analytic rank $r_C \leq 3/2$.
\end{theorem}
\begin{proof}
We adjust \cite{Young:2004} Section 5.5. We use Corollary \ref{cor:1ld} with
\begin{align*}
\phi_\nu(x)&=\left(\frac{\text{sin}(\pi \nu x)}{\pi \nu x} \right)^2 \\
\hat{\phi}_\nu(y) &= \frac{1}{\nu}\left(1- \frac{|y|}{\nu}\right)
\end{align*}
where $\hat{\phi}_\nu$ is supported in $[-\nu,\nu]$ and is the Fourier transform of $\phi_\nu$. Since $\phi_v(0)=1$ and $\phi_v(x)\geq0$, we have $F(0)\geq 1$ and $F(\theta)\geq0$, where we defined $F(\theta)$ in Section \ref{sec:background}. This implies $ \text{ord}_{u=1/q} L(E\otimes\chi_D, u) \leq Z_{\phi_\nu}(\Theta_D)$, so $r_{N,C}\leq \langle  Z_{\phi_\nu}\rangle_{N,C}$. By Corollary \ref{cor:1ld} 
\[ r_{N,C}\leq \int_{O(M)} Z_{\phi_\nu}(\Theta) d\Theta + \mathcal{O}_{E,q}(1/N). \]
By Equation \eqref{eq:rank}, we need to evaluate
\begin{align*}
 \hat{\phi_\nu}(0)+\frac{1}{M}\sum_{n\neq 0} \hat{\phi}_\nu\left(\frac{n}{M}\right) \eta_n&=\frac{1}{\nu}+\frac{2}{M\nu}\sum_{n=1}^{\lfloor M\nu/2\rfloor}\left(1-\frac{2n}{M\nu}\right) \\
 &= \frac{1}{\nu} +\frac{2}{M\nu}\left(\lfloor M\nu/2\rfloor - \frac{2}{M\nu}\sum_{n=1}^{\lfloor M\nu/2\rfloor} n \right) \\
 &=\frac{1}{\nu} +\frac{2}{M\nu}\left(M\nu/2 +\mathcal{O}(1) - \frac{2}{M\nu}\left(\frac{M^2\nu^2}{8}+\mathcal{O}(M\nu)\right) \right)\\
 &=\frac{1}{\nu} + \frac{1}{2} + \mathcal{O}\left(\frac{1}{M\nu}\right).
 \end{align*}
Setting $\nu=1-\epsilon$ for any $\epsilon>0$, we have
\[ r_{N,C} \leq \frac{1}{1-\epsilon} + \frac{1}{2} + \mathcal{O}_{E,q}(1/N)\]
meaning that 
\[ \text{limsup}_{N \rightarrow \infty} r_{N,C} \leq 3/2 + o(1).\]
If the limit exists, then
\[ r_C \leq 3/2.\]
We remark that unlike \cite{Bui&Florea:2016}, we cannot optimize the choice of the test function in order to improve our result as in \cite{Iwaniec&Luo&Sarnak:1999} Appendix A, Corollary 2 since the symmetry is orthogonal.
\end{proof}

Let $\epsilon$ and $\epsilon_D$ be the sign of the functional equation of $L(E,u)$ and $L(E\otimes \chi_D,u)$ respectively. They are both $\pm 1$.
\begin{lemma}\label{lm:r}
The order of the central zero $\emph{ord}_{u=1/q} L(E\otimes\chi_D,u)$ is even or odd depending on whether $\epsilon_D=1$ or $\epsilon_D=-1$ respectively.
\end{lemma}
\begin{proof}
Since $L(E,u)\in\mathbb{Z}[u]$, we also have $L(E\otimes\chi_D,u) \in \mathbb{Z}[u]$ because $\chi_D$ is quadratic. The functional equation is then
\[ L\left(E\otimes\chi_D,\frac{1}{q^2u}\right) = \frac{1}{\epsilon_D u^M q^M} L(E\otimes\chi_D,u)\]
where $M$ is the number of zeros of $L(E\otimes\chi_D,u)$. Applying the functional equation twice shows that $\epsilon_D=\pm1$. \\ \\
If $(1/q)e^{i\theta_j}$ is a zero of $L(E\otimes\chi_D,u)$, then $(1/q)e^{-i\theta_j}$ is also a zero by the functional equation. This means all zeros come in pairs expect those at $u=-1/q$ and $u=1/q$. We have
\[ \epsilon_D = \prod_{i=1}^M -e^{i\theta_j} = (-1)^M (-1)^{\text{ord}_{u=-1/q} L(E\otimes\chi_D,u)}\]
So for example if $\epsilon_D=1$ and $M$ is even, then $\text{ord}_{u=-1/q} L(E\otimes\chi_D,u)$ must be even and $\text{ord}_{u=1/q} L(E\otimes\chi_D,u)$ must be even too. The argument is the same for the other cases.
\end{proof}

\begin{theorem}\label{thm:nv}
At least $12.5\%$ of the family of quadratic twists of an elliptic curve have rank zero and at least $37.5\%$ have rank one as $N\rightarrow \infty$.
\end{theorem}
\begin{proof}
By (1.2) \cite{Bui&Florea&Keating&Roditty-Gershon:2019}
\[ \epsilon_D = \alpha \cdot \epsilon\cdot \chi_{M_E}(D)\]
for $D \in \mathcal{H}^*_N$ where $M_E$ is the product of the primes of multiplicative reduction of $E$, $\alpha=\pm1$ depending only on the degree of $D$, and $\chi_{M_E}$ is the unique quadratic Dirichlet character of conductor $M_E$. We recall that we assume $M_E\neq1$. Since $\chi_{M_E}$ is periodic modulo $M_E$ and since  $M_E\mid N_E$, we have that $\chi_{M_E}$ is constant on $\mathcal{H}_{N,C}$. \\ \\
We recall that we split $\mathcal{H}^*_N=\mathcal{H}_{N,+}\cup\mathcal{H}_{N,-}$ depending on if $\epsilon_D=1$ or $\epsilon_D=-1$ respectively.\\ \\
Let $\chi_{N_E}(D):=\mathds{1}_{(D,N_E)=1}\chi_{M_E}(D)$, which is a Dirichlet character modulo $N_E$ since $M_E\mid N_E$. We have $\chi_{N_E}=\chi_{M_E}$ on $\mathcal{H}^*_N$ since all elements are coprime to $N_E$. By Lemma \ref{lm:size}, $\mathcal{H}_N^*$ becomes equidistributed modulo $N_E$ as $N$ grows,  and since $\chi_{N_E}=1$ on half of $(\mathbb{F}_q[t]/(N_E))^*$  and $\chi_{N_E}=-1$ on the other half, we have
\begin{align}\label{eq:equi}
|\mathcal{H}_{N,+}| \sim |\mathcal{H}_{N,-}| \sim \frac{|\mathcal{H}^*_N|}{2} \quad \text{as }N\rightarrow \infty.
\end{align}
By Theorem \ref{thm:rank}, assuming all limits exist, the average analytic rank is $\leq 3/2$ for both $\mathcal{H}_{N,+}$ and $\mathcal{H}_{N,-}$, since both sets are disjoint unions of congruence classes modulo $N_E$. Then, by Lemma \ref{lm:r}, the rank must be zero for at least $25\%$ of the twists of  $\mathcal{H}_{N,+}$, and for $\mathcal{H}_{N,-}$, the rank must be one for at least $75\%$ of the twists in order to satisfy the bound as $N\rightarrow\infty$. We must divide these quantities by two to conclude.
\end{proof}

\textbf{Remark 1.} 
We have the relation
\[ \epsilon_D = (-1)^M \text{det } \Theta_D=(-1)^{\text{deg}(N_E)} \text{det } \Theta_D\]
where $M$ is the number of zeros of $L(E\otimes\chi_D,u)$. When $\text{deg}(N_E)$ is even, the matrices $\Theta_D$ with $D\in\mathcal{H}_{N,+}$ lie in SO$(M)$ and the matrices $\Theta_D$ with $D\in\mathcal{H}_{N,-}$ lie in $O(M)\backslash SO(M)$. Since $M$ is even, the twists from $\mathcal{H}_{N,+}$ have SO(even) symmetry, and the twists from $\mathcal{H}_{N,-}$ have SO(odd) symmetry. \\ \\
When $\text{deg}(N_E)$ is odd, the matrices $\Theta_D$ with $D\in\mathcal{H}_{N,-}$ lie in SO$(M)$ and the matrices $\Theta_D$ with $D\in\mathcal{H}_{N,+}$ lie in $O(M)\backslash SO(M)$. Since $M$ is odd, the twists from $\mathcal{H}_{N,-}$ have SO(odd) symmetry, and the twists from $\mathcal{H}_{N,+}$ have SO(even) symmetry.  \\ \\
 By \eqref{eq:equi}, the twists from $\mathcal{H}^*_N$ have orthogonal symmetry. See (6) \cite{Young:2004} for a table of statistics of the one-level density of some symmetry types. \\ \\
\textbf{Remark 2.}  The Katz and Sarnak philosophy predicts that Corollary \ref{cor:1ld} should hold without any restriction on the support of $\hat{\phi}$. However, it is important to integrate over the appropriate symmetry group, since it matters when the support of $\hat{\phi}$ is greater than $(-1,1)$. Then, we expect the average analytic rank of the family of twists from $\mathcal{H}^*_N$ to be $1/2$. For the family of twists from $\mathcal{H}_{N,+}$, we expect an average analytic rank of zero, and for $\mathcal{H}_{N,-}$ we expect an average of one. In other words, we expect half of the twists from $\mathcal{H}_N^*$ to have rank zero, and the other half to have rank one.  \\ \\
\textbf{Remark 3.}
The contribution of the primes can be trivially bounded using the Riemann hypothesis (Proposition \ref{prop:rh2})
\[S(N,n):= \sum_{D\in \mathcal{H}^*_N} \sum_{\text{deg}(P)=n} (n/d)(\alpha_P+\overline{\alpha}_P)\chi_D(P)= \mathcal{O}_{E,q}(Nq^{n+N}).\]
The analogue of Hypothesis M from \cite{Fiorilli:2014} is
 \[ S(N,n) = o_{E,q}(q^{n+N}).\]
This bound implies that the contribution of the primes goes to zero as $N\rightarrow \infty$ for any $n$, which is what we need to remove the restriction on the support of $\hat{\phi}$ and to show that the average analytic rank is $1/2$.

\section{Trivial bound for arbitrary order}
We now study the one-level density when the order $\ell \neq 2$ of the Dirichlet characters is coprime to $q$. Using the Lindelöf hypothesis, we obtain the following restriction on the support of the Fourier transform of the test functions
\[ \text{supp } \hat{\phi} \subset \left( -\frac{1}{2},\frac{1}{2}\right).\]
This is the analogue of supp $\hat{\phi} \subset (-1,1)$ for cubic Dirichlet $L$-functions. It is not clear how to use duality to extend the support in this case, as was done for quadratic characters. By restricting to a subfamily (of density zero), David and Güloğlu \cite{chantal:2020} were able to increase the support for cubic Dirichlet $L$-functions over $\mathbb{Q}(\zeta_3)$.
\begin{theorem}\label{thm:ao}
For $\phi \in \mathcal{S}(\mathbb{R})$ an even function such that \emph{supp}$(\hat{\phi})\subset (-1/2,1/2)$, we have for the one-level density of the family of twists of order $\ell\neq 2$ and $(\ell,q)=1$
\[ \langle Z_{\phi}\rangle_{N,C}=\int_{U(M)} Z_{\phi}(\Theta) d\Theta + \mathcal{O}_{E,q,\ell}(1/N)\]
where $M=2N+\emph{deg}(N_E)-4$.
\end{theorem}
\begin{proof}

Let $a>1$ be a divisor of $\ell$ and let $k$ be the smallest integer such that $a \mid q^k-1$, so $k$ depends on $a$. We fix an isomorphism $\psi$ from the $a$-roots of unity in $\mathbb{F}_{q^k}$ to those in $\mathbb{C}$ and we define the $a$-power residue symbol as
\[ \left( \frac{A}{\mathfrak{q} }\right)_{a} := \psi\left( A^{\frac{q^{\text{deg}(\mathfrak{q} )}-1}{a}} \text{ mod }\mathfrak{q}  \right)\]
for $A,\mathfrak{q} \in \mathbb{F}_{q^k}[t]$ where $\mathfrak{q} $ is prime. If $\mathfrak{q} \mid A$, the symbol equals zero by definition. This symbol is a Dirichlet character of modulus $\mathfrak{q}$ of order $a$ over $\mathbb{F}_{q^k}[t]$. Let $\text{Frob}_q$ act on $\mathbb{F}_{q^k}[t]$ by raising each coefficient to the $q$th power. On $\mathbb{F}_{q^k} \subset\mathbb{F}_{q^k}[t]$, this operator permutes the $a$-roots of unity and can be associated via $\psi$ to some $\sigma \in \text{Gal}(\mathbb{Q}(\zeta_a)/\mathbb{Q})$ which generates a normal cyclic subgroup $H\trianglelefteq \text{Gal}(\mathbb{Q}(\zeta_a)/\mathbb{Q})$ of order $k$. Furthermore, we have
\[ \left( \frac{\text{Frob}_q(A)}{\text{Frob}_q(\mathfrak{q} )}  \right)_{a} = \sigma\left(\frac{A}{\mathfrak{q} } \right)_{a}.\]
Now, over $\mathbb{F}_q[t]$, a prime modulus has primitive Dirichlet characters of order $a$ if and only if its degree is divisible by $k$. In other words, the modulus must split totally in the extension $\mathbb{F}_{q^k}[t]$. Moduli of prime powers $\geq 2$ do not have primitive characters of order $a$. For each prime $Q$ such that $k\mid\text{deg}(Q)$, we define $\chi_{a,Q,i}$ for $1\leq i\leq \phi(a)$ to be the collection of characters of order $a$ with conductor $Q$, and they are Gal$(\mathbb{Q}(\zeta_a)/\mathbb{Q})$ conjugates. We extend this definition to composite conductors by multiplicativity. In the following generating series, we consider all primitive characters of order $a \mid \ell$ of conductor of degree $h$ for $a>1$. In the Euler product, we include for each prime all of its possible characters of order $a>1$ dividing $\ell$. Let $|u|=q^{-1/2-\epsilon}$ for any $\epsilon>0$. We have
\begin{align*} \mathcal{L}_P(u)&:=\sum_{h=0}^\infty  \left(\sum_{\substack{\chi \text{ of order } a\mid\ell \\ \chi \text{ primitive} \\ \text{deg}(\text{cond}(\chi))=h \\ (\text{cond}(\chi),N_E)=1 \\ a>1 }}\chi(P)\right)u^h \\
&= \prod_{\substack{Q \text{ prime} \\ Q\nmid N_E}} \left(1+ \sum_{\substack{a\mid\ell \\ a\mid q^{\text{deg}(Q)}-1 \\ a>1}}\sum_{i=1}^{\phi(a)} \chi_{a,Q,i}(P)u^{\text{deg}(Q)}\right) \\
&= \left(\prod_{\substack{Q \text{ prime} }}\prod_{\substack{a\mid\ell \\ a\mid q^{\text{deg}(Q)}-1 \\ a>1}}\prod_{i=1}^{\phi(a)} (1 - \chi_{a,Q,i}(P)u^{\text{deg}(Q)})^{-1}+\mathcal{O}_{\epsilon,\ell,q}(1)\right) \prod_{\substack{Q \text{ prime} \\ Q\mid N_E}} \left(1+ \sum_{\substack{a\mid\ell \\ a\mid q^{\text{deg}(Q)}-1 \\ a>1}}\sum_{i=1}^{\phi(a)} \chi_{a,Q,i}(P)u^{\text{deg}(Q)}\right)^{-1}
\end{align*}
The error comes from the fact that the terms $u^{j \,\text{deg}(Q)}$ with $j\geq 2$  in the product are negligible over $C_\epsilon:|u|= q^{-1/2-\epsilon}$, since there are at most $q^m$ primes of degree $m$ and the sum
\[ \sum_{m=1}^\infty q^m q^{(-1-2\epsilon)m}\]
converges absolutely. We also have
\[ \prod_{\substack{Q \text{ prime} \\ Q\mid N_E}} \left(1+\sum_{\substack{a\mid\ell \\ a\mid q^{\text{deg}(Q)}-1 \\ a>1}}\sum_{i=1}^{\phi(a)} \chi_{a,Q,i}(P) u^{\text{deg}(Q)}\right)^{-1} \ll_{E,\ell,q} 1 \]
over $C_\epsilon$.
Now, we fix $a>1$ dividing $\ell$ and we collect the primes $Q$ having characters of order $a$. Let $k$ be the smallest integer such that $a \mid q^k-1$. When we write ``$\mathfrak{q}$ totally split'' under the product, we mean that the prime in $\mathbb{F}_q[t]$ under $\mathfrak{q}$ splits completely in $\mathbb{F}_{q^k}[t]$. We have
\begin{align*}
\prod_{\substack{Q \text{ prime}  \\ a\mid q^{\text{deg}(Q)}-1}}\prod_{i=1}^{\phi(a)} (1 - \chi_{a,Q,i}(P)u^{\text{deg}(Q)})^{-1}=\prod_{i=1}^{\phi(a)/k} \prod_{\substack{\mathfrak{q} \in \mathbb{F}_{q^k}[t] \\ \mathfrak{q}\text{ prime} \\ \mathfrak{q} \text{ totally split}}} (1-\chi_{a,\mathfrak{q},i}(P)u^{k\,\text{deg}(\mathfrak{q})})^{-1}
\end{align*}
by defining 
\[ \chi_{a,\mathfrak{q},i}(A) := \sigma_i \left( \frac{A}{\mathfrak{q}} \right)_{a} \]
where the $\sigma_i$ are a set of representatives of Gal$(\mathbb{Q}(\zeta_a)/\mathbb{Q})/H$ where $H$ was defined above. Since $\text{Frob}_q(P)=P$, and since a character of conductor $\mathfrak{q}$ of order $a$ over $\mathbb{F}_{q^k}[t]$ restricted to $\mathbb{F}_q[t]$ gives a character of conductor $Q$ of order $a$ over $\mathbb{F}_q[t]$, we have a one-to-one correspondence between these characters. By the $a$-power reciprocity law (\cite{rosen:2013} Theorem 3.5)
\[\prod_{\substack{\mathfrak{q} \in \mathbb{F}_{q^k}[t] \\ \mathfrak{q}\text{ prime} \\ \mathfrak{q} \text{ totally split}}} (1-\chi_{a,\mathfrak{q},i}(P)u^{k\,\text{deg}(\mathfrak{q})})^{-1}=L_{q^k}(\chi_{a,P,i},(-1)^{\frac{q^k-1}{a}\text{deg}(P)}u^k)\prod_{\substack{\mathfrak{q} \in \mathbb{F}_{q^k}[t] \\ \mathfrak{q}\text{ prime} \\ \mathfrak{q} \text{ not totally split}}} (1-\chi_{a,\mathfrak{q},i}(P)u^{k\,\text{deg}(\mathfrak{q})})\]
where the index $q^k$ indicates that the $L$-function is taken over $\mathbb{F}_{q^k}[t]$.
When $\mathfrak{q}$ is not totally split, we have $k \,\text{deg}(\mathfrak{q})\geq 2 \,\text{deg}(Q)$, so the last product converges absolutely on $C_\epsilon$ since the number of primes above a given prime $Q$ is bounded by $\ell$. Finally, by the Lindelöf hypothesis (\cite{David&Florea&Lalin:2019} Lemma 2.6), we have for any $\delta>0$
\begin{align*}
L_{q^k}(\chi_{a,P,i},u^k) &\ll_{q,\delta} q^{k\delta\, \text{deg}(P)} \quad\, \, \, \,\text{for } |u|\leq q^{-1/2}.
\end{align*}
Combining everything gives
\[ \mathcal{L}_P(u) \ll_{q,\delta,\epsilon,E,\ell} q^{\delta\,\text{deg}(P)}\]
on $C_\epsilon$. We now apply Perron's formula
\[ \sum_{\substack{\chi \text{ of order } a\mid\ell \\ \chi \text{ primitive} \\ \text{deg}(\text{cond}(\chi))=h \\ (\text{cond}(\chi),N_E)=1 \\ a>1 }}\chi(P) = \frac{1}{2 \pi i}\oint_{C_\epsilon} \frac{\mathcal{L}_P(u)}{u^{h}} \frac{du}{u}.\]
Since $\mathcal{L}_P(u)$ is entire inside $C_\epsilon$ and since $|1/u^{h+1}| \leq q^{(1/2+\epsilon)h+1}$, we have
\[\sum_{\substack{\chi \text{ of order } a\mid\ell \\ \chi \text{ primitive} \\ \text{deg}(\text{cond}(\chi))=h \\ (\text{cond}(\chi),N_E)=1 \\ a>1 }}\chi(P) \ll_{q,\delta,\epsilon,E,\ell} q^{h/2} q^{\epsilon h} q^{\delta \,\text{deg}(P)} \]
for any $\epsilon,\delta > 0$. Since this holds for any $\ell$, we may isolate the characters of order $\ell$
\[ \sum_{\substack{\chi \text{ of order }\ell \\ \chi \text{ primitive} \\ \text{deg}(\text{cond}(\chi))=h \\ (\text{cond}(\chi),N_E)=1}}\chi(P) \ll_{q,\delta,\epsilon,E,\ell} q^{h/2} q^{\epsilon h} q^{\delta \,\text{deg}(P)}. \]
As in the quadratic case, the average of traces is given by
\[ \langle \text{tr } \Theta^n \rangle_N = \frac{1}{q^n |\mathcal{F}_N|} \sum_{d|n} \sum_{\text{deg}(P)=n/d} (n/d) (\alpha_P^d + \overline{\alpha}_P^d) \sum_{\substack{\chi \text{ of order }\ell \\ \chi \text{ primitive}  \\ \text{deg}(\text{cond}(\chi))=N \\ (\text{cond}(\chi),N_E)=1}}\chi^d(P) \]
where $\mathcal{F}_N$ denotes the family of primitive twists of order $\ell$ of conductor of degree $N$ coprime to $N_E$. Using $|\mathcal{F}_N| \asymp q^N$ and the above, we have
\[ \langle \text{tr } \Theta^n \rangle_N \ll_{q,\epsilon,E,\ell} n^2 q^{n/2} q^{-N/2} q^{\epsilon(n+N)}  \]
since $\chi^2$ isn't principal and higher powers are trivially bounded. This translates to a restriction of $(-1/2,1/2)$ on the support of the test function for the one-level density. The symmetry is unitary since there is no contribution.
\end{proof}
\textbf{Remark.} In this section, the number of zeros of the $L$-functions do not depend on the parity of the characters. However, if we were to remove the elliptic curve, the genus would change depending on the parity of the characters for conductors of fixed degree. See \cite{David&Florea&Lalin:2019} Lemma 2.9 for example. We can control the parity of the characters here by fixing $s:=(\ell,q-1)$ and $\tilde{\chi}$ a character of order $s$ on $\mathbb{F}_q^*$. The indexing of the collection $\{\chi_{a,Q,i}\}_{1\leq i\leq \phi(a)}$ of characters over $\mathbb{F}_q[t]$ described above is now important. We set $\chi_{a,Q,1}$ to be any character from the collection such that it is equal to $\tilde{\chi}^{t_a\text{deg(Q)}}$ when restricted to $\mathbb{F}_q^*$ where $t_a=(\ell,q-1)/(a,q-1)$. Then, we set the other characters as $\chi_{a,Q,i}:=\chi_{a,Q,1}^{b_{a,i}}$ where $b_{a,i}$ is the $i$th element of $(\mathbb{Z}/a\mathbb{Z})^*$. To keep track of the values of the characters over $\mathbb{F}_q^*$, we modify the Euler product as
\[\prod_{\substack{Q \text{ prime} \\ Q\nmid N_E}} \left(1+ \sum_{\substack{a\mid\ell \\ a\mid q^{\text{deg}(Q)}-1 \\ a>1}}\sum_{i=1}^{\phi(a)} \chi_{a,Q,i}(P)\zeta_s^{r t_a b_{a,i}\text{deg}(Q)}u^{\text{deg}(Q)}\right)\]
 where $\zeta_s$ is the $s$th complex root of unity and $0\leq r \leq s-1$. Doing so will change the $L$-function obtained above as
 \[ L_{q^k}(\chi_{a,P,i},(-1)^{\frac{q^k-1}{a}\text{deg}(P)}(\zeta_s^{r t_a b_{a,i}}u)^k)\]
 so it doesn't change the bound we obtained
\[ \sum_{\substack{\chi \text{ of order }\ell \\ \chi \text{ primitive} \\ \text{deg}(\text{cond}(\chi))=h \\ (\text{cond}(\chi),N_E)=1}}\chi(P)\zeta_s^{rb(\chi)} \ll_{q,\delta,\epsilon,E,\ell} q^{h/2} q^{\epsilon h} q^{\delta \,\text{deg}(P)}\]
 where $0\leq b(\chi)\leq s-1$ is the unique number such that
 \[ \chi|_{\mathbb{F}_q^*} = \tilde{\chi}^{b(\chi)}.\]
We can now control the congruence of $b(\chi)$ modulo $s$ by summing over all $r$
\[ \sum_{\substack{\chi \text{ of order }\ell \\ \chi \text{ primitive} \\ \text{deg}(\text{cond}(\chi))=h \\ (\text{cond}(\chi),N_E)=1 \\ b(\chi) \equiv c \text{ mod } s}}\chi(P) = \frac{1}{s}\sum_{r=0}^{s-1} \zeta_s^{-rc}  \sum_{\substack{\chi \text{ of order }\ell \\ \chi \text{ primitive} \\ \text{deg}(\text{cond}(\chi))=h \\ (\text{cond}(\chi),N_E)=1}}\chi(P)\zeta_s^{rb(\chi)}\]
so the bound also holds for this sum. By denoting the genus by $g$ and the degree of the conductor by $h$, we have $g=h-1$ if $\chi$ is odd and $g=h-2$ if $\chi$ is even for Dirichlet $L$-functions. We know that $\chi$ is even if and only if $b(\chi) \equiv 0$ mod $s$, so controlling the congruence gives us control over the genus.
\section{Acknowledgments} 
I would like to thank Chantal David, my doctoral thesis advisor, for the idea of this project, for her support, her help, and the various improvements she suggested. I would also like to thank the Fonds de recherche Nature et technologies du Québec for their financial support.